\documentclass[10pt, journal]{IEEEtran}
\usepackage{graphicx}
\graphicspath{{./imgArxiv/}}
\usepackage{tikz}
\usetikzlibrary{shapes}
\usepackage{pgfplots}
\usepackage{color}
\usepackage{amssymb,amsfonts,amsmath,latexsym,dsfont,bm}

\usepackage{algorithm,algorithmic}

\usepackage{multirow}

\usepackage{subfig}

\usepackage{nccmath}

\usepackage{hyperref}
\hypersetup{
     colorlinks   = true,
     citecolor    = blue,
     linkcolor    = blue,
     urlcolor     = blue
}

\newcommand{\bmo}{\bm{o}}
\newcommand{\bmp}{\bm{p}}
\newcommand{\bmq}{\bm{q}}
\newcommand{\bmu}{\bm{u}}
\newcommand{\bmv}{\bm{v}}
\newcommand{\bmx}{\bm{x}}
\newcommand{\bmy}{\bm{y}}

\newcommand{\bmphi}{\bm{\phi}}
\newcommand{\bmbeta}{\bm{\beta}}

\newcommand{\mcA}{\mathcal{A}}
\newcommand{\mcC}{\mathcal{C}}
\newcommand{\mcI}{\mathcal{I}}
\newcommand{\mcO}{\mathcal{O}}

\begin{document}

\title{Gauss--Newton Optimization for Phase Recovery from the Bispectrum}

% authors
\author{James~L.~Herring, James~Nagy, and~Lars~Ruthotto
\thanks{J.~Herring, Department of Mathematics, University of Houston, Houston, TX 77004, USA email: {\tt herring@math.uh.edu}}%
\thanks{J.~Nagy and L.~Ruthotto, Department of Mathematics, Emory University, Atlanta, GA 30322, email: {\tt \{jnagy,lruthotto\}@emory.edu}}}%

%\markboth{IEEE Transactions on Computational Imaging}%

\maketitle

\begin{abstract}
Phase recovery from the bispectrum is a central problem in speckle interferometry which can be posed as an optimization problem minimizing a weighted nonlinear least-squares objective function. We look at two different formulations of the phase recovery problem from the literature, both of which can be minimized with respect to either the recovered phase or the recovered image. Previously, strategies for solving these formulations have been limited to gradient descent or quasi-Newton methods. This paper explores Gauss--Newton optimization schemes for the problem of phase recovery from the bispectrum. We implement efficient Gauss--Newton optimization schemes for all the formulations. For the two of these formulations which optimize with respect to the recovered image, we also extend to projected Gauss--Newton to enforce element-wise lower and upper bounds on the pixel intensities of the recovered image. We show that our efficient Gauss--Newton schemes result in better image reconstructions with no or limited additional computational cost compared to previously implemented first-order optimization schemes for phase recovery from the bispectrum.  MATLAB implementations of all methods and simulations are made publicly available in the \texttt{BiBox} repository on Github.
\end{abstract}


% Introduction
\section{Introduction}
\label{sec:Intro}

\IEEEPARstart{I}{mage} blurring due to turbulence poses a significant obstacle in many applications. One approach to obtaining high spatial frequency images of an object through turbulent optical systems is speckle interferometry, which is built upon Labeyrie's observation that high spatial frequency information can be recovered from short-exposure images~\cite{Labeyrie1970}. This work provided the means to obtain diffraction-limited reconstructions of an object's Fourier modulus, or power spectrum, using short-exposure, photon-limited data. 

While an object's Fourier modulus may be sufficient in the case of some simple objects, many applications also require recovery of the object's Fourier phase in order to produce high quality images. Thus, phase retrieval is often an essential subproblem when using speckle interferometry. Several methods have been developed to solve the phase retrieval problem, most using relationships defined by high-order statistical correlation measures such as an object's triple correlation or its Fourier counterpart, the bispectrum~\cite{Knox1974, Weigelt1977}. These measures can be collected from the short-exposure data and used to reconstruct an object's phase. 

This work focuses on the problem of phase recovery from an object's collected bispectrum. The broader problem of phase retrieval occurs in numerous engineering fields and the sciences including astronomy, electronmicroscopy, crystallography, and optical imaging. For a recent overview of applications and challenges, we recommend~\cite{Shechtman2015}. The problem of phase recovery from the bispectrum has its origins in astronomical imaging for low-light images in the visible spectrum~\cite{Knox1974, Weigelt1977, Lohmann1983, Lohmann1984, Negrete1996, Wirnitzer1985}. More recently, the problem has been of interest in multiple applications including long-range horizontal and slant path imaging~\cite{Archer2014, Bos2011, Carrano2002} and phase recovery in underwater imaging~\cite{Wen2007, Wen2010}. The general problem of phase retrieval from Fourier measurements also continues to be an active area of research in the signal processing community~\cite{Bendory2017, Bendory2017b, Bendory2018, Bendory2018b, Jaganathan2013, Jaganathan2017, Pinilla2018, Repetti2014, Wang2018}. Many phase retrieval applications considered in the literature are based on the relationship between an object's phase and higher order statistical moments collected from the data (such as the bispectrum) and require solving constrained nonlinear least-squares problems; see, e.g.,~\cite{Bendory2018, Repetti2014, Wang2018}. Thus while we consider only the problem of phase recovery from the bispectrum, the content of this paper is more broadly applicable to the phase retrieval problem in general. 

For phase recovery from the bispectrum, strategies can be separated into two categories: recursive algorithms and weighted least-squares problems. Recursive strategies fix a small set of known phase values near the center of the Fourier domain and use these known phase values along with the collected bispectrum to recursively compute the remainder of the reconstructed phase. Such strategies are well explored in the literature~\cite{Lohmann1983, Lohmann1984, Bos2011, Cho1995, Kang1991,  Matson1991, Meng1990, Northcott1988, Tyler2004}. One limitation of recursive strategies is poor performance when attempting to reconstruct the phase values associated with high frequency information when using noisy data~\cite{Negrete1996}. To improve on this, the phase recovery problem can be reformulated as a weighted least-squares problem minimizing the mismatch between the unknown object phase and the collected phase of the bispectrum. This least-squares problem is typically nonlinear because the bispectrum is collected modulo-$2\pi$ in the range $[-\pi,\pi]$, i.e., the bispectrum is ``wrapped." One can solve this nonlinear least-squares problem~\cite{Negrete1996, Matson1991, Meng1990, Haniff1991, Glindemann1993, Takajo1988, Erdem1992}. Alternatively, techniques have been proposed to unwrap the collected bispectrum, e.g.,~\cite{Goodman1990, Marron1990}, resulting in a linear least-squares problem, but these approaches have been shown to produce inferior results to solving the nonlinear formulation of the problem \cite{Haniff1991}. Previous approaches to solving both the linear and nonlinear least-squares formulations have focused on gradient-based first-order and quasi-Newton methods such as gradient descent and the limited memory Broyden–Fletcher–Goldfarb–Shanno method (L-BFGS)~\cite{Negrete1996, Haniff1991, Glindemann1993}. 

In this paper, we focus on the formulation of phase recovery from the bispectrum as a nonlinear weighted least-squares problem using the wrapped bispectrum. We present Gauss--Newton schemes as an alternative to the gradient-based optimization approaches previously used in the literature. Specifically, we make the following contributions:

\begin{itemize}
\item We implement efficient Gauss--Newton optimization methods for two nonlinear least-squares formulations of the phase recovery problem from the literature. Both of these formulations can be minimized with respect to the recovered phase or recovered image, resulting in four possible formulations. Our implementations exploit sparsity, matrix reorderings, incomplete factorization, and the speed of the fast Fourier transform (FFT) to reduce the cost associated with solving the linear system to calculate the Gauss--Newton step at each iteration of the optimization. The resulting schemes have per iteration costs on the same order of magnitude as gradient descent and quasi-Newton methods like NLCG and L-BFGS but benefit from faster convergence. This results in improved time-to-solution and lower overall computational cost compared with previous approaches. 

\item For the two formulations that are minimized with respect to the resulting image, we also explore the constrained problem with pixel-wise non-negativity constraints on the pixel intensities of the recovered image. We use a projected Gauss--Newton method to enforce these constraints. This strategy improves the quality of the recovered image. Additionally, it eliminates the need for a regularizer enforcing non-negativity and offers the option to use other regularization options while still enforcing non-negativity. We demonstrate this for two common regularizers, a discrete gradient regularizer, and a total variation regularizer. For comparison with constrained gradient-based optimization, we also test projected gradient descent with both regularizers for these two formulations. 

\item We show through numerical experiments that our proposed standard Gauss--Newton and projected Gauss--Newton implementations offer improvements on gradient descent, L-BFGS, and projected gradient descent in both time-to-solution and the quality of the resulting image for the four formulations of the phase recovery problem taken from the literature. We also show that these improvements in the quality of the resulting image are robust to problem parameters including atmospheric turbulence and noise level. 

\item We provide MATLAB implementations for all our methods and for all our numerical simulations in the \texttt{BiBox} repository on Github:
\begin{center}
	\url{https://github.com/herrinj/BiBox}
\end{center}
\end{itemize}

The paper is organized as follows. Section~\ref{sec:PhaseRecovery} gives an overview of the collection of the bispectrum and the formulation of phase recovery from the bispectrum as a weighted nonlinear least-squares optimization problem. Section~\ref{sec:Optimization} details our proposed iterative Gauss--Newton optimization for solving the phase recovery problem from the previous section. Specifically, we introduce expressions for the gradient and Gauss--Newton approximations to the Hessians for four separate formulations of the phase recovery problem and discuss efficient strategies for solving the linear system associated with the Gauss--Newton step within each iteration of the optimization. Lastly, for two of the problem formulations, we extend our implementation to projected Gauss--Newton which allows for phase recovery while simultaneously imposing pixel-wise intensity bounds on the recovered image. Section~\ref{sec:NumericalExperiments} presents numerical experiments. We compare our proposed Gauss--Newton schemes with common first-order optimization methods: gradient descent, projected gradient descent, and quasi-Newton L-BFGS method. We also demonstrate the robustness of our proposed Gauss--Newton schemes for a range of problem parameters. We end with concluding remarks in Section~\ref{sec:Conclusions}.

% Phase Recovery Problem
\section{Phase Recovery Problem}
\label{sec:PhaseRecovery}

Most techniques for recovering an object's phase from speckle image data rely on the object's triple correlation and its Fourier transform, the bispectrum~\cite{Knox1974, Weigelt1977}.  An object's triple correlation is a second-order statistical moment given by measuring an object against two independently shifted copies of itself. For a two dimensional object $o(\bmx)$ with $\bmx \in \mathbb{R}^2$, this is expressed by
\begin{equation*}
o^{TC}(\bmx_1, \bmx_2) = \iint_{-\infty}^\infty  o^*(\bmx) o(\bmx + \bmx_1) o(\bmx + \bmx_2) d\bmx.
\end{equation*}
\noindent Taking the Fourier transform of this and using the convolution property, we get the object's bispectrum, 
\begin{equation*}
O^{(3)}(\bmu, \bmv) = O(\bmu) O(\bmv) O^*(\bmu + \bmv).
\end{equation*}
\noindent Here, $O(\bmu)$ is the Fourier transform of the object and $\bmu, \bmv \in \mathbb{R}^2$ are spatial frequencies. It follows that the phase of the object's bispectrum and the object are related by
\begin{equation}
\label{eq:bispectrumPhase}
\beta(\bmu, \bmv) = \phi(\bmu) + \phi(\bmv) - \phi(\bmu + \bmv)
\end{equation}
\noindent where $\phi(\bmu)$ is the phase of $O(\bmu)$ and $\beta(\bmu,\bmv)$ is the phase of the object's bispectrum corresponding to the triplet $(\bmu,\bmv, \bmu + \bmv)$. This expression provides a deterministic relationship that provides the basis for various algorithms for phase recovery from the bispectrum. 

We highlight several characteristics of the relationship in~(\ref{eq:bispectrumPhase}). First, the object's bispectrum is unknown in practice and is estimated using the average data bispectrum from a series of short-exposure images of the object, $i_k(\bmx)$ for $k = 1,\ldots,N$. To recover the object's Fourier modulus, or power spectrum, an additional set of short-exposure images of an appropriate reference star, $s_k(\bmx)$ for $k = 1,\ldots,N$, is also necessary; see, e.g.,~\cite{Negrete1996}. Typically, the data bispectrum is collected modulo-$2\pi$ or ``wrapped." This necessitates accounting for the modulus within the phase recovery problem or unwrapping the bispectrum before phase recovery; see,  e.g.,~\cite{Goodman1990, Marron1990}. We opt for the first strategy. 

The size of the data bispectrum depends on the number of discrete $(\bmu, \bmv, \bmu + \bmv)$ triplets for which the bispectrum is collected. Computational considerations and the optics of the problem make it inadvisable to range $\bmu$ and $\bmv$ over the entire Fourier plane. In practice, the number of discrete coordinates $\bmu$ and $\bmv$ used to determine these triplets is restricted by two radii: $\bmu$ indices are restricted within a larger, ``recovery radius," and $\bmv$ indices are restricted within a smaller radius~\cite{Buscher1993, Nakajima1988, Roggemann1992}. The value of these radii vary depending on seeing conditions and parameters; we discuss the values we select in the discussion on data setup in Sec.~\ref{sec:NumericalExperiments}. For computational purposes, the indices for the set of discrete $(\bmu, \bmv, \bmu + \bmv)$ triplets are computed once and stored in an indexing structure which vectorizes the accumulation of the data bispectrum. This structure can be implemented to exploit the symmetries in the phase for real-valued images in Fourier space, which increases computational efficiency. We base our indexing structure on the work presented by~\cite{Tyler2004}.

%%% Phase Recovery Schemes %%%
\subsection*{Phase Recovery Schemes}

Numerous algorithms for phase recovery from the \mbox{modulo-$2\pi$} bispectrum have been proposed which build upon the relationship in (\ref{eq:bispectrumPhase}). These can be separated into two categories: recursive algorithms; see, e.g.,~\cite{Lohmann1983, Lohmann1984, Bos2011, Cho1995, Kang1991,  Matson1991, Meng1990, Northcott1988, Tyler2004}, and nonlinear least-squares formulations; see, e.g.,~\cite{Negrete1996,  Matson1991, Meng1990, Haniff1991,Glindemann1993, Takajo1988}. We look at the second category. Specifically, we look at four nonlinear least-squares formulations proposed in the literature.

We first establish some notation. Let $\bmphi_{\rm true} \in \mathbb{R}^n$ be the vector of true object phase values we aim to recover. Here, the dimension $n$ is determined by seeing conditions and signal-to-noise (SNR) considerations. Let $\bmbeta_{\rm true} \in \mathbb{R}^m$ be the corresponding true, unwrapped phase of the data bispectrum collected for $m$ distinct $(\bmu,\bmv,\bmu+\bmv)$ triplets. These two vectors are related by
\begin{equation}
\label{eq:Axb}
\bmbeta_{\rm true} = A \bmphi_{\rm true}
\end{equation}
\noindent where $A \in \mathbb{R}^{m \times n}$ is a sparse matrix with three non-zeros per row: two 1's and one -1 corresponding to the signs of phase elements in (\ref{eq:bispectrumPhase}). The indexing structure used to accumulate the data bispectrum also contains the information used to construct this matrix. The relationship in~(\ref{eq:Axb}) provides the basis of the algorithms for recovering the unknown object phase, $\bmphi \in \mathbb{R}^n$, from the phase of the collected data bispectrum, $\bmbeta \in \mathbb{R}^m$. The SNR for each entry in $\bmbeta$ varies, so we define a diagonal weighting matrix $W \in \mathbb{R}^{m \times m}$ with positive diagonal entries determined by the SNR of the bispectrum phase elements~\cite{Negrete1996}. Let $W^{1/2}$ be the diagonal matrix with the square roots of the SNR weights on the diagonal.

Using these definitions, we introduce several optimization formulations for recovering an object's phase from the phase of the data bispectrum from the literature. All of the formulations involve fitting the data through the minimization of an appropriate nonlinear least-squares objective function. Here, we introduce two such objective functions from the literature. The first fits the object's phase to the collected data bispectrum by solving the following minimization problem~\cite{Haniff1991}:   
\begin{equation}
\label{eq:Haniff}
\min_{\phi} \left\{E_1(\bmphi)  = \frac{1}{2} \| W^{1/2} \operatorname{mod}_{2\pi} \left(\bmbeta - A\bmphi \right) \|^2_2\right\}.
\end{equation}  
\noindent Here, the modulus is introduced because the phase of the data bispectrum is collected modulo-$2\pi$. This introduces several considerations. It causes $E_1(\bmphi)$ to be non-convex with periodic local minima every $2\pi$. Additionally, $E_1(\bmphi)$ is periodically discontinuous with jumps every $2\pi$ where the misfit wraps from $0$ to $2\pi$. 

An alternative formulation to~(\ref{eq:Haniff}) exists which avoids the issue of non-differentiability~\cite{Haniff1991}. It exploits the identity $e^{i\theta} = \cos(\theta) + i \sin(\theta)$ to avoid the modulus. The resulting optimization problem is
\begin{equation}
\label{eq:Haniff2}
\begin{split}
\min_{\bmphi} \bigg \{ E_2(\bmphi) = &\frac{1}{2} \| W^{1/2}\left(\cos \bmbeta - \cos A \bmphi \right)\|_2^2 \\
&+  \frac{1}{2}\|W^{1/2}\left(\sin \bmbeta - \sin A \bmphi\right) \|_2^2 \bigg \} , 
\end{split}
\end{equation}

Both $E_1(\bmphi)$ and $E_2(\bmphi)$ are minimized with respect to the unknown object phase, $\bmphi$. One shortcoming is that these objective functions solely minimize the data misfit between the recovered phase and the collected data bispectrum, i.e, they are blind to how the recovered phase impacts the recovered object when combined with the object's power spectrum. An alternative is to optimize with respect to the resulting image, $\bmo$, i.e., to reformulate (\ref{eq:Haniff}) and (\ref{eq:Haniff2}) as $E_1(\bmphi(\bmo)) = E_1(\bmo)$ and $E_2(\bmphi(\bmo)) = E_2(\bmo)$, respectively. This idea was originally proposed for $E_1(\bmo)$ in~\cite{Glindemann1993}. Extending the idea to $E_2(\bmo)$ is new to our knowledge.

Optimizing with respect to the resulting image is attractive for a number of reasons. In many applications, the image of the recovered object is the ultimate goal which necessitates solving the phase recovery problem. It then makes intuitive sense that the optimization should take into account this recovered image. Additionally, since the objective functions for both~(\ref{eq:Haniff}) and~(\ref{eq:Haniff2}) are highly non-convex, optimizing with respect to the resulting image may allow us to converge to different minima, potentially resulting in improved recovery of the object. Lastly, optimizing with respect to the recovered image allows us to impose desirable characteristics on the solution image using appropriate regularization and bound-constraints. Such characteristics may include non-negativity or smoothness of the recovered object. With this in mind, we incorporate regularization and bound constraints into optimization formulations for $E_1(\bmo)$ and $E_2(\bmo)$ and consider the problems
\begin{equation}
\label{eq:Glindemann}
\min_{\bmo \in \mcC} \left\{ E_1(\bmo) + \alpha R(\bmo) \right\}
\end{equation}
\noindent and 
\begin{equation}
\label{eq:Glindemann2}
\min_{\bmo \in \mcC} \left\{ E_2(\bmo) + \alpha R(\bmo) \right\}\,.
\end{equation}
Here, $\mcC$ denotes a closed, convex set enforcing element-wise bound constraints on the recovered object $\bmo$, for example, non-negativity. The operator $R(\bmo)$ is a regularizer subject to the weighting parameter $\alpha > 0$ for introducing desirable characteristics in the solution such as smoothness, sparsity, or non-negativity. In Sec.~\ref{sec:NumericalExperiments}, we look at three potential regularizers. For the unconstrained problem ($\mcC = \mathbb{R}^n$), we explore a previously proposed penalty regularizer encouraging non-negative solutions~\cite{Glindemann1993}. We also solve the problem with non-negativity constraints, which we enforce using projected Gauss--Newton. For this formulation, we implement and test two regularizers: a quadratic regularizer with a discrete gradient operator, $\nabla_h$, and a nonlinear total variation regularizer popular in many imaging applications~\cite{Rudin1992}. Other potential options include the identity operator, $p$-norm based regularizers, or hybrid regularization techniques.

%%% Optimization %%%
\section{Optimization}
\label{sec:Optimization}

We now look at methods for minimizing the four optimization problems introduced in Sec.~\ref{sec:PhaseRecovery}.
Previous work has relied on first-order and quasi-Newton optimization methods such as gradient descent, nonlinear conjugate gradient (NLCG), and the limited memory implementation of the Broyden--Fletcher--Goldfarb--Shanno (L-BFGS) method~\cite{Negrete1996, Haniff1991, Glindemann1993}. For further references on these methods, see~\cite{NocedalWright1999}. Our work implements efficient, Gauss--Newton optimization schemes. To begin our discussion, we provide a brief overview of the Gauss--Newton method. After this, we discuss extending the method to projected Gauss--Newton as a way to enforce element-wise bound constraints on the solution. Lastly, we narrow our focus and discuss specific approaches for solving the linear system to determine the Gauss--Newton step for each of the objective functions proposed in the previous section. Solving this linear system is the most computationally intensive step in Gauss--Newton optimization and represents the most significant contribution of our work. For (\ref{eq:Haniff}) and (\ref{eq:Haniff2}), we introduce an efficient factorization based approach which uses sparse matrix reordering and factorization to solve the linear system directly. For (\ref{eq:Glindemann}) and (\ref{eq:Glindemann2}), we solve for the Gauss--Newton step iteratively. We also extend the optimization for (\ref{eq:Glindemann}) and (\ref{eq:Glindemann2}) to include regularization.

%%% Gauss--Newton %%%
\subsection*{Gauss--Newton Method}

The Gauss--Newton method is a common optimization scheme for solving unconstrained nonlinear least-squares problems such as~(\ref{eq:Haniff}) and~(\ref{eq:Haniff2})~\cite{NocedalWright1999, Beck2014}. Starting from an initial guess, the method iteratively updates the computed solution with a step towards minimizing the objective function. At each iteration, both the direction and length of this step must be computed. 

The direction of this step at each iteration is determined by solving a linear system. For an arbitrary objective function $E(\bmy)$ at the current iterate $\bmy$, this equation is 
\begin{equation}
\label{eq:GNstep}
 H_E(\bmy) \bmp = - \nabla_{\bmy} E(\bmy)\,.
\end{equation}
\noindent Here, $\nabla_{\bmy} E(\bmy)$ is the gradient, and $H_E(\bmy) \approx \nabla^2_{\bmy} E(\bmy)$ is a symmetric positive semi-definite approximation to the Hessian; see, e.g., p.~246 in \cite{NocedalWright1999}. Both the gradient and the approximation to the Hessian must be updated at each iteration of the optimization. The system~(\ref{eq:GNstep}) can be solved using direct methods or an appropriate iterative solver, and solving this system is the most significant cost in the Gauss-Newton method. Our work focuses on implementing efficient strategies for solving this system, and we use both direct and iterative solvers for~(\ref{eq:GNstep}) depending on which of the four formulations of the phase recovery problem we consider.

After choosing the step direction, we choose the step length via a backtracking Armijo line search~\cite{NocedalWright1999}. The Armijo line search selects a step length $0  < \eta \leq 1$ by iteratively backtracking from the full Newton step $\eta = 1$ until a step length is found which guarantees a sufficient reduction of the objective function. Sufficient reduction is determined by the Armijo condition, 
\begin{equation*}
E(\bmy + \eta \bmp) \leq E(\bmy) + c \eta \nabla_{\bmy} E(\bmy)^\top \bmp
\end{equation*}
\noindent where $c \in (0,1)$ is a reduction constant. We use $c = 1 \times 10^{-4}$ as recommended in~\cite{NocedalWright1999}. Note that if the data residual is small in the neighborhood of the solution, the Gauss--Newton approximation to the Hessian is more accurate, and we expect to take the full step $\eta = 1$. This reduces the number of backtracking line search iterations required per outer Gauss--Newton iteration and lowers the cost of the optimization.

The convergence of the Gauss-Newton method depends on the nonlinearity of the objective function and the norm of the data residual. In the best case, it can be shown to have near-quadratic convergence in the neighborhood of a minimizer. This makes it attractive compared to the linear convergence of gradient-based methods and the super-linear convergence of quasi-Newton methods like L-BFGS. We note that the backtracking Armijo line search we use does not satisfy the Wolfe conditions necessary for rigorous convergence results using Gauss--Newton, but this does not pose a problem in practice as observed in Sec.~\ref{sec:NumericalExperiments}. For a more in-depth discussion of the convergence properties of Gauss--Newton and the other optimization methods; see, e.g.,~\cite{NocedalWright1999}.

%%% Projected Gauss--Newton %%%
\subsection*{Projected Gauss--Newton Method}

The standard Gauss--Newton method is sufficient for the unconstrained optimization problems in (\ref{eq:Haniff}) and (\ref{eq:Haniff2}). However, additional considerations are required to solve constrained problems. For the formulations (\ref{eq:Glindemann}) and (\ref{eq:Glindemann2}), we want to enforce element-wise bound constraints on the solution object. To do this, we use the projected Gauss--Newton method~\cite{Haber2014}. Projected Gauss--Newton incorporates constraints by combining Gauss--Newton with projected gradient descent. The goal of this combination is to benefit from the strengths of both methods: fast convergence in the neighborhood of the solution for Gauss--Newton and straightforward implementation of bound constraints for projected gradient descent. 

We implement projected Gauss--Newton as follows. Given a feasible initial guess, the method separates the set of optimization variables into two sets at each iteration. The active set, $\mcA$, contains the indices of the variables for which the bound constraints are active, i.e., variables on the boundary of the feasible region. The inactive set, $\mcI$, contains the indices of the variables for which the bound constraints are inactive, i.e., variables on the interior of the feasible region. Let the vectors $\bmx_{\mcA}$ and $\bmx_{\mcI}$ be vectors containing the entries of $\bmx$ at the indices defined by the active and inactive sets, respectively. 

On the inactive set, we take a Gauss--Newton step $\bmp_{\mcI}$ computed by solving a projected version of (\ref{eq:GNstep}). This can be written
\begin{equation*}
\left( I_{\mcI} H_E(\bmy) I_{\mcI} \right) \bmp_{\mcI} = - I_{\mcI} \nabla_{\bmy} E(\bmy),
\end{equation*}
\noindent where $I_{\mcI}$ is a projection operator. It is computed by modifying an identity matrix by setting diagonal entries in the columns corresponding to indices in $\mcA$ to $0$. This projected linear system can be solved using either a direct method or an iterative method with an appropriate preconditioner.

On the active set, we take a projected gradient descent step given by
\begin{equation*}
\bmp_{\mcA} = - I_{\mcA} \nabla^P_{\bmy} E(\bmy)\,.
\end{equation*}
Here, $I_{\mcA} = I - I_{\mcI}$ and $\nabla^P_{\bmy} E(\bmy)$ is the projected gradient where the projection sets to zero entries of the gradient that would cause the updated solution to violate the bound constraints. 

We then combine the steps on the inactive and active sets to get the full projected Gauss--Newton step, $\bmp$, given by
\begin{equation*}
\bmp = \bmp_{\mcI} + \gamma \bmp_{\mcA}.
\end{equation*}
\noindent Here, the parameter $\gamma > 0$ is introduced to reconcile the difference in scale between the two steps. Following the recommendation in~\cite{Haber2014}, we set $\gamma = \frac{\|\bmp_{\mcI}\|_{\infty}}{\| \bmp_{\mcA} \|_{\infty}}$. 

Lastly, we determine the length of the step using a projected Armijo line search. It follows the same backtracking strategy presented previously for standard Gauss--Newton but uses a modified Armijo condition,
\begin{equation*}
E(Q(\bmy + \eta \bmp)) \leq E(\bmy) + c \eta \nabla^P_{\bmy} E(\bmy)^\top \bmp,
\end{equation*}
\noindent where $Q$ is a projection onto the feasible region and $\nabla^P_{\bmy} E(\bmy)$ is the projected gradient. Note that due to the projected gradient descent step on the active set, we expect projected Gauss--Newton to require more line search iterations than standard Gauss--Newton. This is particularly relevant when many bound-constraints are active, and we often observe this in practice. Also, variables may enter and leave the active and inactive sets with each update to the solution, so the active and inactive sets must be updated at each iteration of the optimization.

The convergence behavior of projected Gauss--Newton depends on the objective function and the number of variables in the active and inactive sets. For the best case where all variables lie in the inactive set, the method reverts to standard Gauss--Newton and has the corresponding convergence properties. In the worst case scenario where all variables are in the active set, convergence behavior is that of gradient descent, i.e., linear convergence. Our numerical experiments in Sec.~\ref{sec:NumericalExperiments} show that for the phase recovery problems in this paper, we observe more of the former.

%%% Gauss--Newton Step %%%
\subsection*{Efficient Approaches for computing the Gauss--Newton Step}
We now consider solving the linear system~(\ref{eq:GNstep}) associated with the Gauss--Newton step for the four formulations of the problem. We begin by considering the optimization problem associated with  $E_1(\bmphi)$ in equation~(\ref{eq:Haniff}). The gradient and the Gauss--Newton approximation to the Hessian are given by
\begin{equation}
\label{eq:HaniffDerivatives}
\begin{split}
\nabla_{\bmphi} E_1(\bmphi) &= -A^\top W \operatorname{mod}_{2\pi} \left(\bmbeta - A \bmphi \right)\\
H_{E_1(\bmphi)} &= A^\top W A.
\end{split}
\end{equation}
\noindent For both derivatives, we ignore the modulus during differentiation. This strategy is used in the literature and proves effective in practice, but it may impact the optimization~\cite{Haniff1991}. We also note that the approximation to the Hessian, $H_{E_1(\bmphi)}$, is independent of the variable $\bmphi$ and remains constant for each Gauss--Newton iteration.
Furthermore, the matrix is large, sparse, and structured. It has a large number of zero rows and columns due to SNR considerations and symmetries in the phase for real-valued images. The object's phase can only be recovered at points for which data, i.e., the phase of the bispectrum, has been collected. This is restricted with a problem specific `recovery radius' determined by seeing conditions and the optics system~\cite{Buscher1993, Nakajima1988, Roggemann1992}. Thus, columns and rows corresponding to entries outside this radius are set to zero. Furthermore, the phase for real-valued images displays a symmetry through D.C. That is, if we let D.C. be the origin, $(0,0)$, then the phase at point $(i,j)$ is the negative of the phase at point $(-i,-j)$, i.e., $-\phi(i,j) = \phi(-i,-j)$. Exploiting this allows us to further increase sparsity, and makes sparse, direct methods an attractive option for solving the linear system. Importantly, any factorization can be computed once offline and then reused for each iteration of the optimization. 

\begin{figure*}[!t]
\centering
\subfloat{\includegraphics[width = 0.8\textwidth]{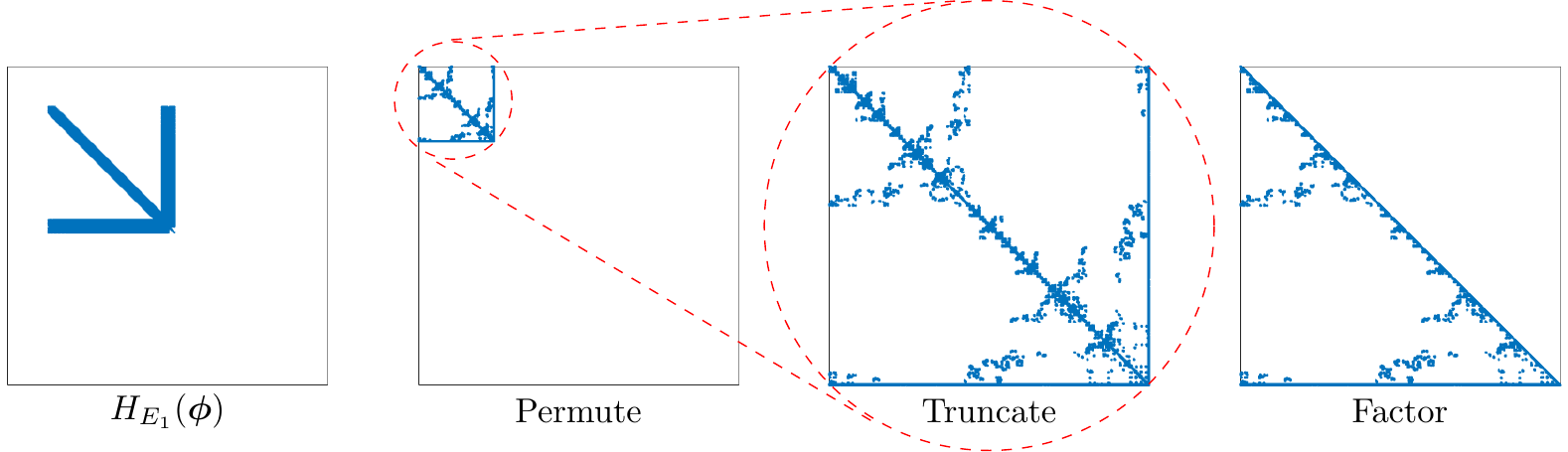}}
\caption{This figure shows the permutation and factorization strategy used for the approximate Hessians of $E_1(\bmphi)$ and $E_2(\bmphi)$. From \textit{left} to \textit{right}, we start with the full Hessian, permute using a symmetric approximate minimum degree factorization, truncate out the zero rows and columns, and perform an incomplete Cholesky factorization on the truncated matrix. We use the resulting factors and the permutation to solve for the (\ref{eq:GNstep}) at each iteration of the Gauss--Newton optimization.}  
\label{fig:cholesky}
\end{figure*}

With that in mind, we implement the following strategy. Prior to the Gauss--Newton optimization, we compute and factor the approximate Hessian offline. We perform a symmetric approximate minimum degree permutation to shift the non-zero rows and columns to the upper left-hand corner of the matrix~\cite{Davis2004}. This permutation is a heuristic designed to minimize the fill-in when the permuted matrix is factored. Then, we extract the $\tilde{n} \times \tilde{n}$ sub-matrix corresponding to the non-zero columns and rows of the permuted matrix. Here, the size $\tilde{n}$ corresponds to the number of phase-values to be recovered and is small due to the SNR considerations and symmetries previously mentioned, $\tilde{n} < n$. The extracted sub-matrix is also symmetric positive definite. On this sub-matrix, we perform a zero fill-in incomplete Cholesky factorization and store the factors~\cite{Manteuffel1980,Saad2003}. The symmetric approximate minimum degree permutation helps to minimize the amount of information lost with the zero fill-in strategy. A full Cholesky factorization is also possible, but it may lead to some loss of sparsity due to fill-in and does not improve results for the problems presented in this paper. The stored incomplete Cholesky factors can then be used to solve for the Gauss--Newton step at each iteration of the standard Gauss--Newton method. Computing the factorization once offline reduces the cost of solving~(\ref{eq:GNstep}) at each iteration. The cost of solving the linear system associated with the Gauss--Newton step is thereby reduced to solving two sparse, triangular systems for a cost of $\mcO (\tilde{n}^2)$. Note that in practice, $\mcO (\tilde{n}^2)$ is a pessimistic bound on the cost of inverting the triangular factors because the incomplete Cholesky factorization preserves their sparsity. Figure~\ref{fig:cholesky} gives a visual step-by-step of the strategy.  

In the case of $E_2(\bmphi)$, the gradient and the Gauss--Newton approximation to the Hessian are
\begin{equation}
\label{eq:Haniff2Derivatives}
\begin{split}
\nabla_{\bmphi} E_2(\bmphi) &= A^\top W D_1(\bmphi) \\
H_{E_2}(\bmphi) &= A^\top W D_2(\bmphi) A\,.
\end{split}
\end{equation}
The matrices $D_1(\bmphi)$ and $D_2(\bmphi)$ are diagonal matrices with diagonal entries given by 
\begin{equation*}
\begin{split}
\texttt{diag}(D_1(\bmphi)) &= \cos \bmbeta \odot \sin A \bmphi - \sin \bmbeta \odot \cos A\bmphi \\
\texttt{diag}(D_2(\bmphi)) &= \cos \bmbeta \odot \cos A \bmphi + \sin \bmbeta \odot \sin A \bmphi,
\end{split}
\end{equation*}
\noindent where $\odot$ denotes the Hadamard product or component-wise multiplication of two vectors. Unlike~(\ref{eq:HaniffDerivatives}), the Gauss--Newton approximation to the Hessian here is not independent of the phase due to the diagonal matrix $D_2(\bmphi)$, and makes solving the system in (\ref{eq:Haniff2Derivatives}) more computationally expensive. However, we note that if $\operatorname{mod}_{2\pi}(\bmbeta - A \bmphi) = 0$, i.e., we perfectly match the recovered phase to the phase of the bispectrum, then the diagonal entries of $D_2(\bmphi)$ are $1$, i.e, $D_2(\bmphi)$ becomes the identity. In practice, this is not true for noisy data and large residuals, but it serves as motivation for omitting the diagonal matrix $D_2(\bmphi)$ from the approximate Hessian. This makes $H_{E_1}(\bmphi) = H_{E_2}(\bmphi)$, and we can use the same permutation and factorization strategy as we use for (\ref{eq:Haniff}) when computing the Gauss--Newton step for (\ref{eq:Haniff2}). This is effective in practice, as shown in the numerical experiments. We also remark that although we use the same approximate Hessian for both formulations, the objective functions and gradients differ resulting in different solutions. 
 
For the optimization problems based on $E_1(\bmo)$ and $E_2(\bmo)$, that is solving (\ref{eq:Glindemann}) and (\ref{eq:Glindemann2})  respectively,
it is necessary to differentiate with respect to the resulting image. Expressions for the gradients and the Gauss--Newton approximations to the Hessians follow from (\ref{eq:HaniffDerivatives}) and (\ref{eq:Haniff2Derivatives}) via the chain rule with additional gradients and derivatives for the regularization term. For the objective function in (\ref{eq:Glindemann}), we get 
\begin{equation}
\label{eq:GlindemannDerivatives}
\begin{split}
\nabla_{\bmo} E_1(\bmphi) &= -\frac{\partial \bmphi^*}{\partial \bmo}A^\top W \operatorname{mod}_{2\pi} \left(\bmbeta - A \bmphi \right) + \alpha \nabla_{\bmo} R(\bmo) \\
H_{E_1}(\bmphi) &= \frac{\partial \bmphi^*}{\partial \bmo} A^\top W A \frac{\partial \bmphi}{\partial \bmo} + \alpha  \nabla_{\bmo}^2 R(\bmo)
\end{split}
\end{equation}
\noindent where $\frac{\partial \bmphi}{\partial \bmo}$ is the derivative of the phase with respect to the resulting object and $\frac{\partial \bmphi^*}{\partial \bmo}$ is its adjoint. Similarly, the derivatives for the objective function in (\ref{eq:Glindemann2}) are given by 
\begin{equation}
\label{eq:Glindemann2Derivatives}
\begin{split}
\nabla_{\bmo} E_2(\bmphi) &= \frac{\partial \bmphi^*}{\partial \bmo} A^\top W D_1(\bmphi) + \alpha \nabla_{\bmo} R(\bmo)\\
H_ {E_2}(\bmphi) &= \frac{\partial \bmphi^*}{\partial \bmo} A^\top W D_2(\bmphi) A \frac{\partial \bmphi}{\partial \bmo} + \alpha  \nabla_{\bmo}^2 R(\bmo).
\end{split}
\end{equation}
To derive expressions for $\frac{\partial \bmphi}{\partial \bmo}$ and its adjoint, we observe that the object's phase can be expressed as a function of the object by
\begin{equation*}
\bmphi(\bmo) = \arctan \left( \frac{\operatorname{Im}(\mathcal{F}\bmo)}{\operatorname{Re}(\mathcal{F}\bmo)} \right), 
\end{equation*} 
\noindent where $\mathcal{F}$ is a 2D Fourier transform matrix. Differentiating this with respect $\bmo$, we get an expression for the action of the adjoint operator in the direction $\bmq$, 
\begin{equation*}
\begin{split}
\frac{\partial \bmphi^*}{\partial \bmo} (\bmq) &= \frac{\operatorname{Re}(\mathcal{F}\bmo) \odot \operatorname{Im}(\mathcal{F}\bmq) - \operatorname{Im}(\mathcal{F}\bmo)\odot  \operatorname{Re}(\mathcal{F}\bmq)}{|\mathcal{F}\bmo|^2} \\
&= \frac{\operatorname{Im} (\overline{\mathcal{F} \bmo} \odot \mathcal{F}\bmq)}{\overline{\mathcal{F} \bmo} \odot \mathcal{F} \bmo}\\
&= \operatorname{Im} \left( \frac{\mathcal{F} \bmq}{ \mathcal{F} \bmo} \right).
\end{split}
\end{equation*} 
\noindent Here, the Hadamard product, division, square, and conjugation are all taken component-wise. This expression is also given in~\cite{Glindemann1993}. 

To evaluate the Gauss--Newton approximations to the Hessians, it is necessary to compute the forward operator, $\frac{\partial \bmphi}{\partial \bmo}$. This operator is not included in previous work using first-order optimization methods. It can be computed and applied using a matrix-free approach. In the direction $\bmq$, it is given by
\begin{equation*}
\begin{split}
\frac{\partial \bmphi}{\partial \bmo} (\bmq) &=  \medmath{\operatorname{Im} \left( \mathcal{F} \left(\frac{\operatorname{Re}(\mathcal{F} \bmo) \odot \bmq }{|\mathcal{F}\bmo|^2} \right)\right) -  \operatorname{Re} \left( \mathcal{F} \left(\frac{\operatorname{Im}(\mathcal{F} \bmo) \odot \bmq}{|\mathcal{F}\bmo|^2} \right)\right)} \\
&= \operatorname{Im} \left( \mathcal{F} \left( \frac{\overline{\mathcal{F} \bmo} \odot \bmq}{\overline{\mathcal{F} \bmo} \odot \mathcal{F}\bmo} \right)\right) \\
&= \operatorname{Im} \left( \mathcal{F} \left( \frac{\bmq}{\mathcal{F} \bmo} \right)\right).
\end{split}
\end{equation*}
Again, all operations are taken component-wise. In practice, some additional considerations are necessary for both the forward and adjoint operators. Storing the 2D FFT matrices is inefficient and infeasible, so the operators are passed as function handles and evaluated as matrix-vector products in a specific direction, $\bmq$. To avoid division by zero, we set the gradient equal to zero when the denominator for indices where $\mathcal{F} \bmo$ is equal to zero. It is also important to scale the Fourier transformations appropriately. Lastly, we note that our code includes an adjoint test to verify the agreement of the forward and adjoint operators. 

Sparse, direct methods are not a viable option when solving for the Gauss--Newton step using~(\ref{eq:GlindemannDerivatives}) and~(\ref{eq:Glindemann2Derivatives}) because the 2D Fourier transform matrices in $\frac{\partial \bmphi}{\partial \bmo}$ and its adjoint are dense, complex matrices and should not be formed explicitly. Instead, we pass the approximate Hessians for $E_1(\bmo)$ and $E_2(\bmo)$ as function handles that evaluate the matrices' product with a vector. Both matrices are symmetric, so separate implementations of their transposes are unnecessary. Many iterative methods that require only matrix-vector products can then be used to solve~(\ref{eq:GNstep}). We use the conjugate gradient method \cite{Bjorck1996, golub2013matrix}. We minimize the cost by solving the system to low accuracy, typically $1\times 10^{-1}$, following the example of~\cite{NocedalWright1999}. This makes solving (\ref{eq:GNstep}) quite efficient. Potentially, we could increase efficiency further by introducing an appropriate preconditioner for the conjugate gradient method. This represents future work.

%%% Numerical Experiments %%%
\section{Numerical Experiments}
\label{sec:NumericalExperiments}

We now compare our proposed methods with methods found in the literature for phase recovery from the bispectrum in several experiments. First, we solve (\ref{eq:Haniff}) and (\ref{eq:Haniff2}) to recover an object's phase and compare results for our proposed Gauss--Newton strategy with gradient descent and L-BFGS. For these formulations of the problem, no regularization is included. We then run experiments with (\ref{eq:Glindemann}) and (\ref{eq:Glindemann2}), optimizing with respect to the recovered object. Here, we compare gradient descent, projected gradient descent, L-BFGS, standard Gauss--Newton, and projected Gauss--Newton. For gradient descent, L-BFGS, and standard Gauss--Newton, we use the penalty term regularizer defined by Eq.~12 in~\cite{Glindemann1993} to encourage non-negativity in the solution. For the constrained problem with strict pixel-wise non-negativity in the solution image, we use projected Gauss--Newton to enforce the bound constraints and compare this to projected gradient descent, a common gradient-based approach for solving the constrained problem. For both projected methods, we test two regularization options for $R(\bmo)$: a quadratic regularizer using a discrete gradient operator and a total variation regularizer. One could also consider a bound constrained approach to L-BFGS~\cite{Byrd1995}, but we do not explore that in this work.

We organize the experiments in the following way. We begin by discussing the setup of the phase recovery problem including simulation of the speckle data and an initial guess. This is followed by a discussion of method parameters for the optimization including line search parameters, stopping criteria, and regularization parameter selection. Then, we show the results for phase recovery when minimizing (\ref{eq:Haniff}), (\ref{eq:Haniff2}), (\ref{eq:Glindemann}), and (\ref{eq:Glindemann2}). These results are split into two sets of experiments. The first set compares all the various optimization schemes for each of the four problem formulations for a problem with fixed parameters. This shows the utility of the proposed Gauss--Newton approach. The second set shows the robustness of the proposed Gauss--Newton approach by testing it across a range of problem parameters.

We compare the results of the various experiments by looking at several values. When comparing optimization methods, we compare the relative change in the objective function with respect to the initial guess (ROF). For an objective function~$E(\cdot)$ at the iterate~$\bmy^{(k)}$ with initial guess~$\bmy^{(0)}$, this is given by $\frac{E(\bmy^{(k)})}{E(\bmy^{(0)})}$. We also compare the relative error of the recovered object (RE), number of optimization iterations (Its.), total optimization time (Time), and CPU time per optimization iteration (Time/It.). Lastly, we track the average number of line search iterations (inner iterations) per optimization iterate (outer iterations) during the optimization (LS/It.). The abbreviations listed in parentheses are used to reference these values in tables and figures. 

We use the following notation to refer to each method in tables and figures: GD for gradient descent, PGD for projected gradient descent, L-BFGS for the limited memory Broyden--Fletcher--Goldfarb--Shanno method, GN for standard Gauss--Newton, and PGN for projected Gauss--Newton. When regularization is used for (\ref{eq:Glindemann}) and (\ref{eq:Glindemann2}), we add a suffix of $+$, $\nabla_h$, or $TV$ in tables and figures to denote the penalty, discrete gradient, and total variation regularizers, respectively. For example, GD$+$ denotes gradient descent with the penalty regularizer.

%%% Data Setup %%%
\subsection*{Data Setup}

To setup test problems, we simulate speckle imaging data for a known true object, $\bmo_{\rm true}$. We use a $256 \times 256$ image of a satellite as the true object, which can be obtained from \cite{nagy2004iterative} (this test image is widely in the literature to evaluate algorithms for image restoration problems; see, e.g., \cite{RoWe96}). We then generate $100$ frames of short exposure data of the object and the reference star for a chosen Fried parameter, $D/r_0$. Higher values of this parameter correspond to more atmospheric turbulence and blurrier images. The approach we use, which is implemented in the MATLAB package {\tt IR\,Tools} \cite{IRtools2018} (software can be obtained from \url{https://github.com/jnagy1/IRtools}), is described in detail in \cite{MatlabWFS2010}. Each data set is generated using a different seed for the random number generator to guarantee independence of the randomness of the two data sets. The object data is scaled to include $3 \times 10^6$ photo-events per data frame, and zero-mean Gaussian noise with a chosen standard deviation $\sigma_{rn}$ is added. The reference star data is scaled to $5000$ photo-events per frame. From this data, we recover the object's power spectrum and accumulate the data bispectrum at a set of $(\bmu, \bmv, \bmu + \bmv)$ triplets. The set of triplets is collected using selected `recovery radius', denoted by $R$, and a smaller radius of $5$, i.e., the indices satisfy $|\bmu| <= R$, $|\bmu+\bmv| <= R$, and $|\bmv| \leq 5$. These triplets are stored in the indexing structure described in Sec.~\ref{sec:PhaseRecovery} and are used to generate the matrix $A$. The choices for the Fried parameter $D/r_0$, recovery radius $R$, and standard deviation of the Gaussian noise $\sigma_{rn}$ vary for different numerical experiments, so values for these parameters are listed in the subsections corresponding to the relevant set of experiments. All other parameters are fixed as listed here for all numerical experiments.

For all four formulations of the phase recovery problem, (\ref{eq:Haniff}), (\ref{eq:Haniff2}), (\ref{eq:Glindemann}), and(\ref{eq:Glindemann2}), the optimization methods are sensitive to the choice of initial guess for the object phase or image. To generate an appropriate initial guess for (\ref{eq:Haniff}) and (\ref{eq:Haniff2}), we use the phase recovered by a single iteration of the recursive algorithm from \cite{Tyler2004}. We denote this by $\bmphi_{init}$. For (\ref{eq:Glindemann}) and (\ref{eq:Glindemann2}) in the unconstrained case, this is combined with the recovered power spectrum of the data to produce an initial guess for the recovered object, $\bmo_{init}$. The relative error for this initial guess varies due to the random generation of the speckle data but is typically $0.75-0.85$ depending on the problem parameters. When using projected Gauss--Newton, the initial guess for the recovered objected must be within the feasible region, i.e., it must have strictly non-negative pixel intensities. To obtain a non-negative guess for the object, we project $\bmo_{init}$ using a method which finds the nearest non-negative image such that the sum of the pixel intensities is identical to $\bmo_{init}$; this is referred to as an energy preserving constraint, and the implementation we use can be found in the {\tt IR\,Tools} package \cite{IRtools2018}. We denote this projected initial guess by $\hat{\bmo}_{init}$. In practice, we also `bump' this projected guess off the bound constraint by adding $\epsilon = 1 \times 10^{-4}$ to ensure the active set is empty for the first iteration of the projected Gauss--Newton method, i.e., no pixels intensities are exactly zero to start the optimization. We note that this projection changes both the power spectrum and phase of the initial guess. Another option for the projection is a single iteration of the error-reduction algorithm \cite{Fienup1982, Negrete1996b}. Images for the true object, average speckle data frame, initial guess from the recursive algorithm, and projected initial guess can be seen in Fig.~\ref{fig:images}.

%%% Parameter Selection %%%
\subsection*{Parameter selection}

The gradient descent, projected gradient descent, L-BFGS, Gauss--Newton, and projected Gauss--Newton methods require setting several parameters. These include stopping criteria for the methods and parameters for the Armijo line search. Additionally, formulations (\ref{eq:Glindemann}) and (\ref{eq:Glindemann2}) depend on the regularization parameter used for the objective functions. We discuss selections for these parameters here.

For stopping criteria, we monitor three values. The first two are the change in the objective function value and the $2$-norm of the difference between successive iterates. The method stops when these two values fall below tolerances of $1\times 10^{-4}$ and $1 \times 10^{-4}$, respectively. Additionally for the Gauss--Newton and projected Gauss--Newton methods, we use an approximation to the Newton decrement as a stopping criteria~\cite{Boyd2004}. The true Newton decrement monitors the reduction in the curvature of the Hessian near a minimizer and is defined by
\begin{equation*}
( \nabla_{\bmy} E(\bmy)^\top \nabla_{\bmy}^2 E(\bmy)^{-1} \nabla_{\bmy} E(\bmy) )^{1/2}\,.
\end{equation*}
We approximate substituting our Gauss--Newton approximation to the Hessian, $H(\bmy)$, instead of $\nabla_{\bmy}^2 E(\bmy)$. Note  that Gauss--Newton step update is given by $\bmp = - H(\bmy)^{-1} \nabla_{\bmy} E(\bmy)$, so we can cheaply compute the approximate Newton decrement by $(-\nabla_{\bmy} E(\bmy)^\top \bmp )^{1/2}$. Both Gauss--Newton and projected Gauss--Newton stop if either the Newton decrement falls below a tolerance of $1\times 10^{-3}$ or if the above criteria for the objective function and difference between successive iterates are met. 

The line search strategies differ for gradient descent, projected gradient descent, and L-BFGS compared to Gauss--Newton and projected Gauss--Newton. The full Newton step is perfectly scaled, and it follows that $\eta = 1$ is a logical initial choice for the step length for Gauss--Newton and projected Gauss--Newton at each iteration when the line search is called. In contrast, the step updates for gradient descent, projected gradient descent, and L-BFGS do not have natural scaling. To account for this, we use an adaptive line search strategy. It works as follows. If the line search succeeds in its first iteration, i.e., it takes the longest step possible, then the line search doubles its initial guess for the step length at the next optimization iterate, i.e. the next time the line search is called. If the line search backtracks, i.e., it takes a shorter step than its initial guess, then it uses the accepted, shorter step as its initial guess for the step length the next time it is called. This adaptive strategy prevents the gradient-based optimization strategies from repeatedly taking initial line search steps which are too large or too small. This is important because each line search iteration requires reevaluating the objective function, which is costly. The adaptive strategy used for gradient descent, projected gradient descent, and L-BFGS helps to reduce the number of line search iterations per optimization iteration (and therefore reduce cost). This provides a fair grounds for comparison with Gauss--Newton and projected Gauss--Newton methods.

For (\ref{eq:Glindemann}) and (\ref{eq:Glindemann2}), we include regularization on the recovered object. This requires selecting the regularization parameter $\alpha$. Choosing this parameter is a challenging subproblem and an active area of research. Popular methods for selecting it include the unbiased predictive risk method, L-curve, and generalized cross-validation (GCV)\cite{Hansen1998, Vogel2002}. For our simulated problem where the true object is known, we solve both (\ref{eq:Glindemann}) and (\ref{eq:Glindemann2}) for a wide range of parameters and choose the parameter which minimizes the relative error of the recovered image. For the unconstrained problem where the penalty-term regularizer is used, this results in $\alpha = 1 \times 10^3$ for both (\ref{eq:Glindemann}) and (\ref{eq:Glindemann2}). For the constrained problem with the discrete gradient regularizer, $\alpha = 1 \times 10^{-2}$ provides the best solutions for both objective functions, and for the constrained problem with the total variation regularizer, we use  $\alpha = 1 \times 10^{4}$ . We note that these regularization parameters are sensitive to the data parameters and problem set up. Also, for real-world data where the true solution is unknown, this method for regularization parameter selection is infeasible. A more sophisticated strategy for regularization parameter selection represents future work, but by using the best case scenarios for choosing regularization parameters, we are able to provide a fair comparison between the computational efficiency of the various methods.

%%% Comparison of optimization schemes %%%
\subsection{Comparison of Gauss--Newton and gradient-based optimization}
\label{ssec:Comparison}

\begin{figure*}[!t]
\centering
\subfloat{\includegraphics{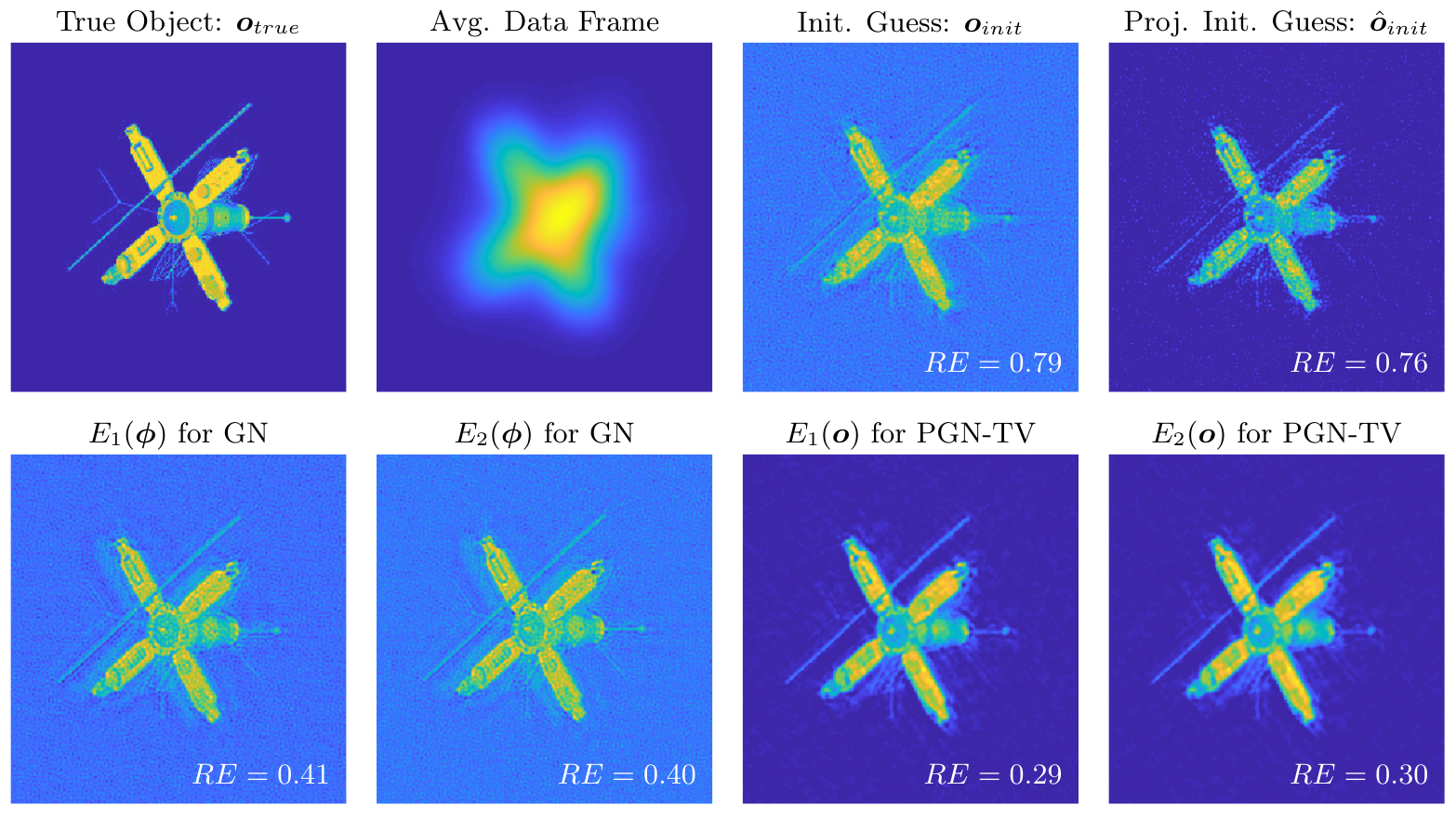}}
\caption{The top row shows the true object, data, and initial guesses used to set up the phase recovery problem. The bottom row shows the best recovered images of the satellite. For $E_1(\bmphi)$ and $E_2(\bmphi)$, the results are shown for standard Gauss--Newton with no regularization on the solution phase. For $E_1(\bmo)$ and $E_2(\bmo)$, the results are shown for projected Gauss--Newton using a total variation regularizer. The relative error of the initial guesses and solution images are printed in each image's bottom right-hand corner.} 
\label{fig:images}
\end{figure*}

\begin{figure*}[!t]
\centering
\subfloat{\includegraphics{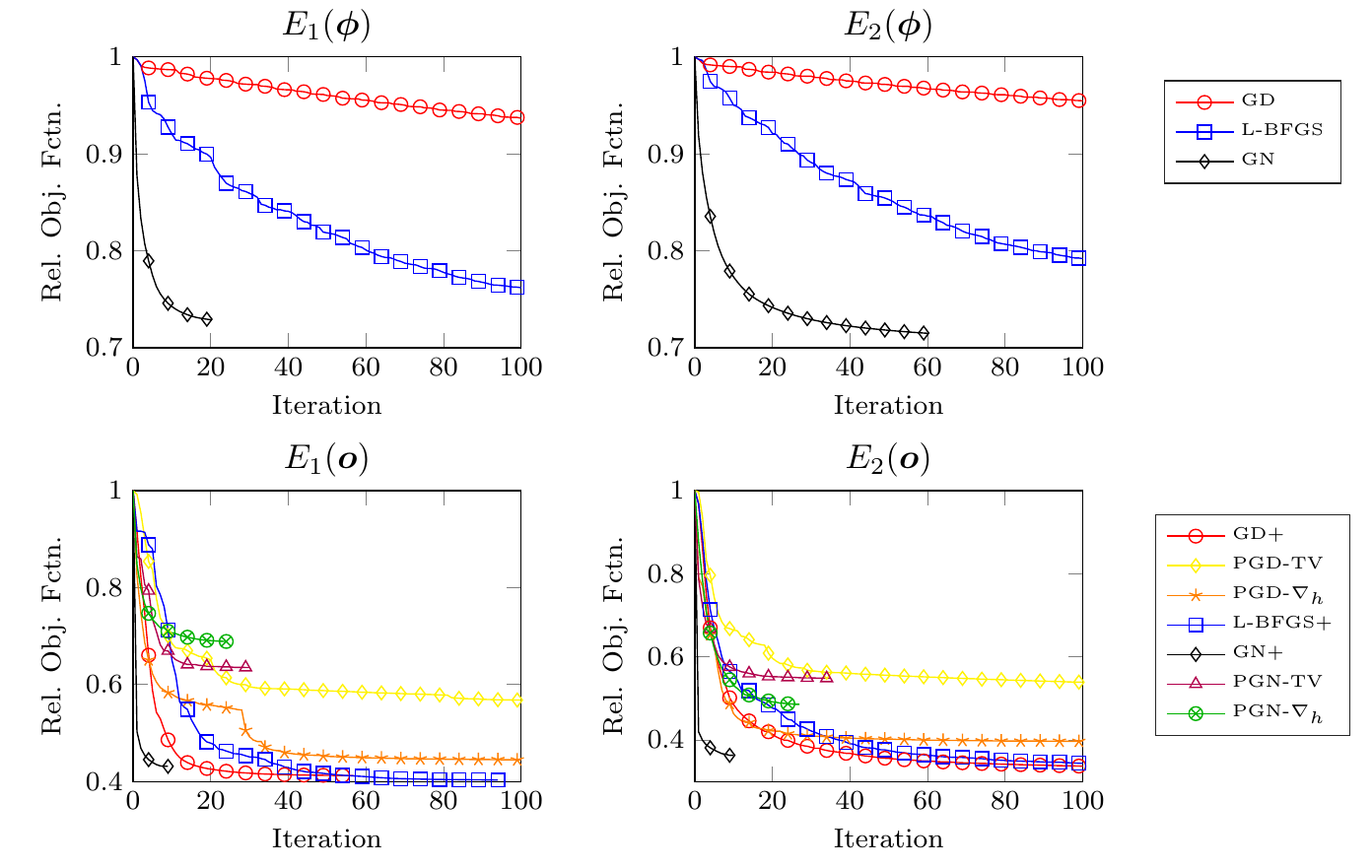}}
\caption{This figure shows convergence plots of the relative objective function for all four formulations of the problem. Note that the Gauss--Newton approaches converge to their minima in fewer iterations than the gradient-based approaches. Also, note that when minimizing $E_1(\bmo)$ and $E_2(\bmo)$, projected gradient descent and projected Gauss--Newton converge to different minima than the other approaches due to the different regularizers, different regularization parameters, and constraints.}
\label{fig:ROF}
\end{figure*}

\begin{table}[!t]
\renewcommand{\arraystretch}{1.3}
\caption{Optimization results for the four formulations of the phase recovery problem. From \textit{left} to \textit{right}, the columns show the optimization method used, minimum relative objective function (ROF), minimum relative error (RE), number of iterations (Its.), CPU time (Time), CPU time per iteration (Time/It.), and line search iterations per optimization iteration (LS/It.) All values were averaged over 50 separate simulated problems.}
\label{tab:results}
\begin{tabular}{c}
\includegraphics{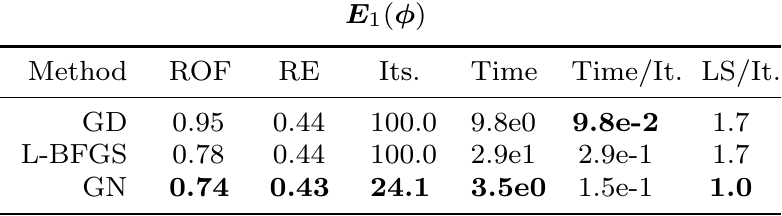} \\
\includegraphics{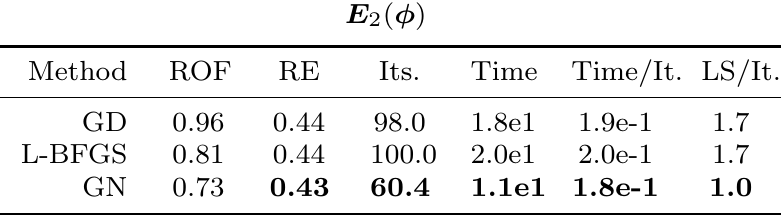} \\ 
\includegraphics{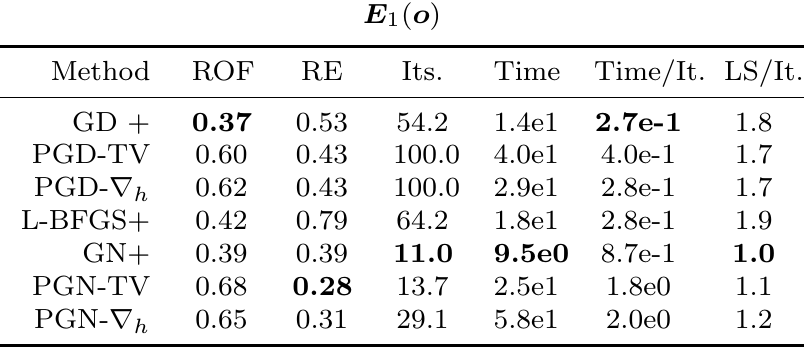} \\ 
\includegraphics{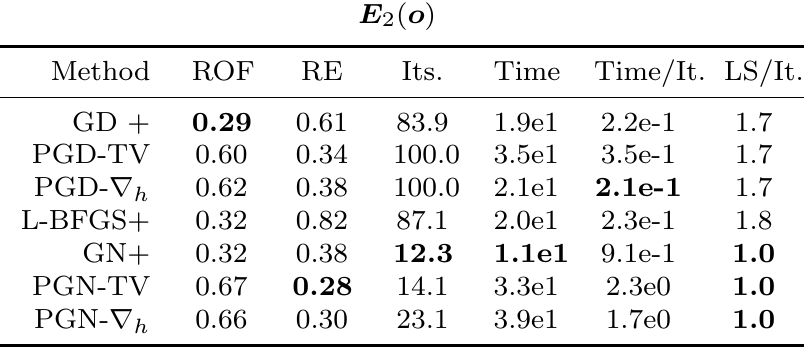} 
\end{tabular}
\end{table}

To compare our proposed Gauss--Newton schemes with previous gradient-based optimization approaches, we solved $50$ different phase recovery problems with independently generated speckle data using all four formulations of the phase recovery problem. For these problems, we used the data parameters and problem setup outlined previously and set the Fried parameter $D/r_0 = 30$, recovery radius $R = 96$, and standard deviation of the Gaussian noise $\sigma_{rn} = 5$. Table~\ref{tab:results} provides results detailing the cost of the optimization and quality of the resulting solutions averaged over all $50$ problems. Fig.~\ref{fig:ROF} compares the convergence of the various optimization schemes using the relative reduction in the objective function for a single problem. Also, the best recovered images for the single problem are displayed in Fig.~\ref{fig:images}. To comment on the results, we separate the discussion for the formulations of the problem which solve for the phase, (\ref{eq:Haniff}) and (\ref{eq:Haniff2}), and those which solve for the resulting object, (\ref{eq:Glindemann}) and (\ref{eq:Glindemann2}).

Looking at the results for $E_1(\bmphi)$ and $E_2(\bmphi)$ in Table~\ref{tab:results}, we see that the Gauss--Newton method outperforms both gradient descent and L-BFGS in terms of minimizing the objective function. This results in marginally lower relative error values for the solution object for the Gauss--Newton method, and all methods offer marginal improvements on the initial guess provided by the recursive algorithm. This is true for both formulations of the problem. In terms of cost, the gradient-based methods require less time per iteration due to the additional computational cost of solving for the step direction at each iteration of the Gauss--Newton method. However, this additional cost is offset by the faster convergence of the Gauss--Newton method (it requires significantly fewer iterations for convergence than either gradient-based optimization scheme). This is observable in Fig.~\ref{fig:ROF}. Overall, the time-to-solution for~(\ref{eq:Haniff}) and~(\ref{eq:Haniff2}) using Gauss--Newton is reduced for both formulations, and as mentioned above, this speed-up is coupled with better solutions in terms of relative error. This suggests that Gauss--Newton optimization is a preferable option to the gradient-based methods for these formulations of the problem. 

The results for (\ref{eq:Glindemann}) and (\ref{eq:Glindemann2}) differ from the two formulations above which consider only the recovered phase. Both standard Gauss--Newton with the penalty term regularizer and projected Gauss--Newton with the total variation and discrete gradient regularizers offer significant improvements in the quality of the recovered object. This can be seen from the minimum relative error values in Table~\ref{tab:results}, where both Gauss--Newton approaches offer significant improvements on the relative error of the initial guess. Projected gradient descent also reduces the relative error of the recovered object for both regularization options, although not as much as the projected Gauss--Newton approach. However, the table shows that gradient descent and L-BFGS  with the penalty regularizer actually worsen the relative error of the recovered image compared to the initial guess despite reducing the objective function. This emphasizes the importance of appropriate regularization and bound constraints in driving the solution to an appropriate minimum. The recovered images also show the advantages of the Gauss--Newton approaches, where the regularization terms in (\ref{eq:Glindemann}) and (\ref{eq:Glindemann2}) reduce ringing and graininess in the solution images in Fig.~\ref{fig:images} compared to the initial guess. 

When comparing the different optimization methods for solving (\ref{eq:Glindemann}) and (\ref{eq:Glindemann2}), we note that the different regularizers and regularization parameters used make the relative objective function a suboptimal way of comparing the methods, see Fig.~\ref{fig:ROF}. Instead, we note that projected Gauss--Newton with non-negativity constraints and the total variation regularizer achieves the best solution images in terms of relative error out of all the methods. However, this improvement comes at a cost as the projected Gauss--Newton approach is the slowest of the four methods. This represents a trade-off, but the method may be preferable if quality of the recovered object is more important than time-to-solution. Also, the projected Gauss--Newton approach allows for the use of various regularizers, so an appropriate regularizer can be chosen for a specific application to optimize solution quality. 

Using standard Gauss--Newton to solve (\ref{eq:Glindemann}) and (\ref{eq:Glindemann2}) with the penalty term regularizer offers the best compromise of speed and quality for our problem. It improves on solution quality and time-to-solution compared to all the gradient-based methods. It is also faster than the projected Gauss--Newton approach, but the resulting images have slightly higher relative error. 

Overall, we observe that for all four formulations of the problem, the standard and projected Gauss--Newton optimization schemes offer improvements in recovered solution images over the gradient descent, L-BFGS, and projected gradient descent methods. Comparing the formulations of the problem, (\ref{eq:Glindemann}) and (\ref{eq:Glindemann2}) are preferable to (\ref{eq:Haniff}) and (\ref{eq:Haniff2}) due to additional regularization and constraints on the solution image that these formulations offer. This results in improvements in the relative error of the recovered image, and the solution images in Fig.~\ref{fig:images} are less grainy and have less ringing. In terms of cost, the faster convergence for standard Gauss--Newton schemes results in lower time-to-solution for all formulations. This is not the case when solving (\ref{eq:Glindemann}) and (\ref{eq:Glindemann2}) with constraints using projected Gauss--Newton. However, the introduction of constraints and appropriate regularization using this approach gives the best-quality images in terms of relative error. This may make this strategy attractive despite the increased cost.

%%% Robustness %%%
\subsection{Robustness of Gauss--Newton schemes}
\label{ssec:Robustness}
\begin{table}[!t]
\renewcommand{\arraystretch}{1.3}
\caption{Optimization results showing the average relative error of the recovered object using the proposed Gauss--Newton schemes for various values of Fried parameter $D/r_0$, recovery radius $R$, and standard deviation of the Gaussian noise $\sigma_{rn}$. From \textit{left} to \textit{right}, the columns show the value of the problem parameter being tested and the relative errors for the initial guess guess, projected initial guess, solution to $E_1(\bmphi)$ using GN, solution to $E_2(\bmphi)$ using GN, solution to $E_1(\bmo)$ minimized using PGN-TV, and the solution to $E_2(\bmo)$ minimized using PGN-TV. All relative error values are averaged over $10$ separate problems.}
\label{tab:results2}
\begin{tabular}{c}
\includegraphics{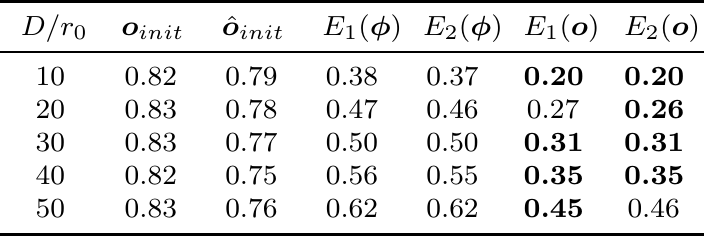} \\
\includegraphics{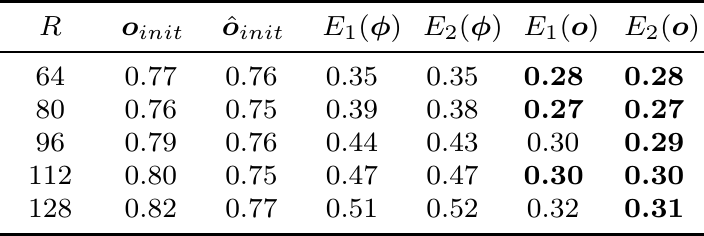} \\ 
\includegraphics{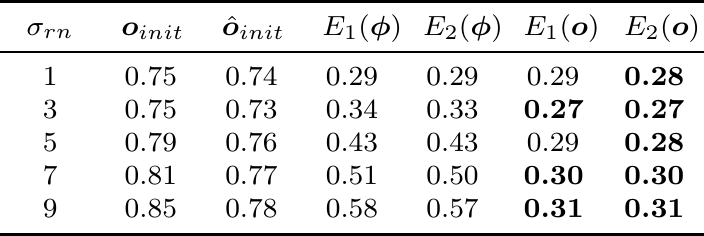} 
\end{tabular}
\end{table}
\begin{figure*}[!t]
\centering
\subfloat{\includegraphics{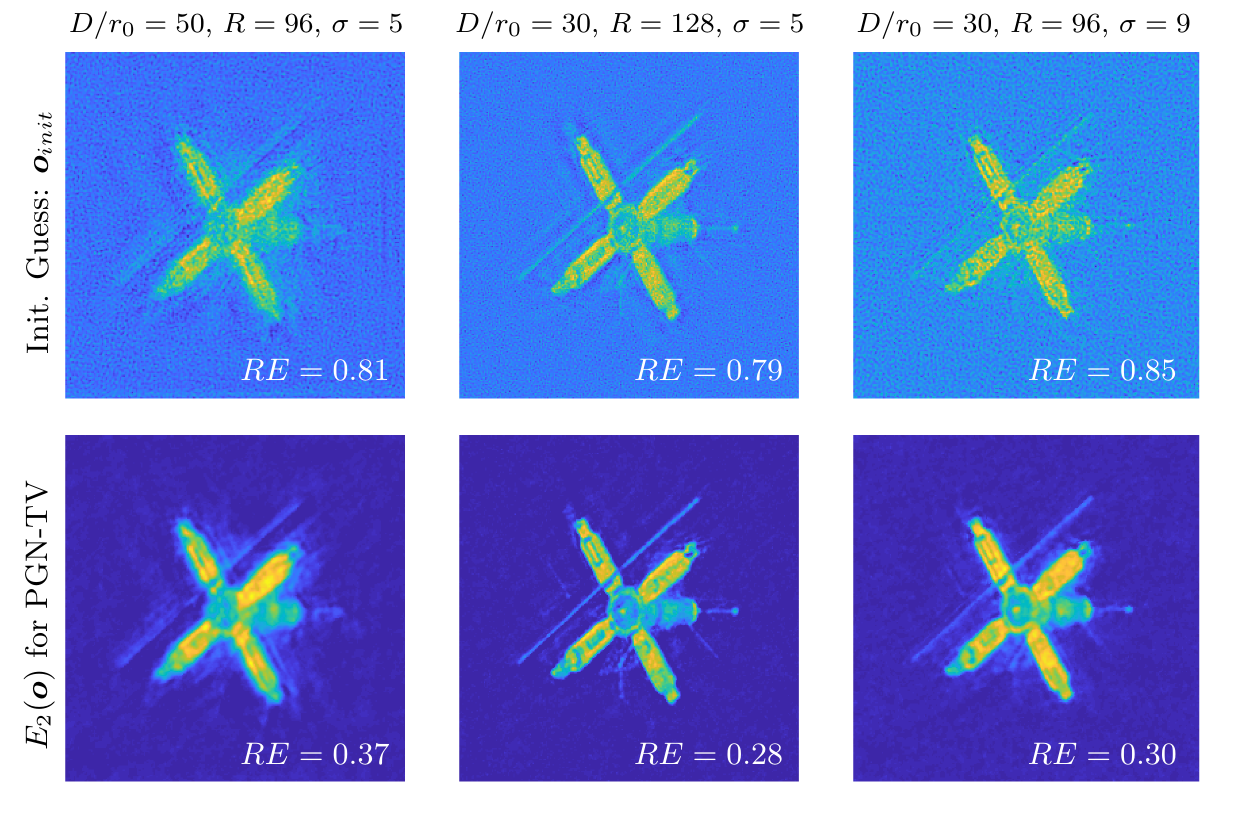}}
\caption{The rows from \textit{top} to \textit{bottom} show the initial guesses given by the recursive algorithm and optimization solutions for $E_2(\bm{o})$ using projected Gauss--Newton with the total variation regularizer (PGN-TV). The columns from \textit{left} to \textit{right} indicate the method parameters used for each problem. The relative error of each image is printed in its bottom right-hand corner.}\label{fig:images2}
\end{figure*}
We also tested the robustness of the our proposed Gauss--Newton schemes by solving several phase recovery problems over a range of parameter values for the Fried parameter $D/r_0$, recovery radius $R$, and standard deviation of the Gaussian noise $\sigma_{rn}$. We set up the problems in the following way. To test the performance with respect to atmospheric blur, we fixed the recovery radius $R = 96$ and standard deviation $\sigma_{rn} = 5$ and solved the problem for Fried parameters $D/r_0 = 10, 20, 30, 40,$ and~$50$. Here, increasing values for $D/r_0$ indicate more atmospheric turbulence and blurrier data. For each value of $D/r_0$, we averaged results over $10$ problems. We then tested the methods for a range of recovery radii, $R = 64, 80, 96, 112,$ and $128$ while fixing $D/r_0 = 30$ and $\sigma_{rn} = 5$. Increasing the recovery radius introduces values of the data bispectrum and recovered phase corresponding to higher spatial frequencies. This presents a tradeoff as data for higher spatial frequencies can potentially enable better recovery of the object being imaged but is more negatively affected by data noise. As with the Fried parameter, we solved $10$ problems for each different value of $R$. Lastly, we solved the phase recovery over a range of values for the standard deviation of the Gaussian noise, $\sigma_{rn} = 1,3,5,7,$ and~$9$ for fixed $D/r_0 = 30$ and $R = 96$. This corresponds to increased noise in the problem data. As with the previous setups, we solved $10$ problems for each value of $\sigma_{rn}$.  

We solved the problems for each of the four objective functions using the respective Gauss--Newton approach that performed best in terms of relative error in the previous set of experiments. For (\ref{eq:Haniff}) and (\ref{eq:Haniff2}), this was the Gauss--Newton approach using the direct permutation and factorization approach described in \ref{sec:Optimization}, denoted GN in Table~\ref{tab:results}. For (\ref{eq:Glindemann}) and (\ref{eq:Glindemann2}), we used projected Gauss--Newton with with the total variation regularizer, denoted PGN-TV. 

We compare the results using the average relative error of the recovered object from solving each objective function across the range problem parameters. This is contrasted with the average relative error of the initial guess, $\bmo_{init}$, from the recursive algorithm and its projection, $\hat{\bmo}_{init}$. These relative error values are displayed in Table~\ref{tab:results2} with the best relative error results highlighted. The performance of the optimization (reduction of the relative objective function, iterations, and time-to-solution, etc.) is similar to the previous set of experiments, so we omit these values. We display the initial guesses and solution images minimizing (\ref{eq:Glindemann2}) using projected Gauss--Newton with the total variation regularizer for three separate problems in Fig.~\ref{fig:images2}. As expected, the results show that the average relative error of the recovered objects increases as the problems become more difficult, i.e., blurrier or noisier data.  Also, phase recovery using the nonlinear least-squares formulations results in better solutions than the initial guesses based on the recursive algorithm. This is likely due to the recursive algorithm's limited ability to recover high spatial frequency information and is evident over the range of problem parameters. The results also reinforce the conclusion that the non-negativity constraints and regularization options available for (\ref{eq:Glindemann}) and (\ref{eq:Glindemann2}) make these formulations of the phase recovery problem preferable to (\ref{eq:Haniff}) and (\ref{eq:Haniff2}), especially if the quality of the recovered object is paramount. This is well illustrated by comparing the initial guesses and solution images using projected Gauss--Newton with the total variation regularizer in Fig.~\ref{fig:images2}.

%%% Conclusions %%%
\section{Conclusions}
\label{sec:Conclusions}

This paper revisits multiple formulations for phase recovery from the bispectrum, a central problem in speckle interferometry. The formulations considered lead to weighted nonlinear least-squares problems which can be solved for either the phase or the object itself. Previous approaches in the literature focused on gradient-based optimization schemes for solving these least-squares problem, including gradient descent and L-BFGS. 

In this work, we implement efficient Gauss--Newton schemes to solve the phase recovery problem. We implement these schemes for four formulations of the problem taken from the literature. To reduce the computational cost of solving the linear system associated with the Gauss--Newton step, we develop tailored approaches for each formulation which exploit the structure and sparsity of the problem. For two formulations of the problem which solve for the solution image, we also extend standard Gauss--Newton to projected Gauss--Newton to allow for element-wise non-negativity constraints on the solution object, a desirable characteristic in many imaging applications. Additionally, the non-negativity constraints within the optimization allows flexibility to introduce additional regularization on the problem. We show this by running numerical experiments with both total variation and discrete gradient regularizers, but the formulation allows for other options. 

Our numerical experiments show that our implementations offer improvements in terms of time-to-solution and quality of solution compared to previously used gradient-based approaches. Our Gauss--Newton schemes produce recovered objects with lower relative error and reduce image artifacts in the resulting images. Furthermore, our implementations achieve these results in less time than the previously implemented first-order methods in most cases. In the cases where Gauss--Newton optimization is more expensive than gradient-based approaches, the increased cost is offset by significantly better solutions. We also show that the improvements in solution quality are robust across a range of problem parameters, particularly for formulations~(\ref{eq:Glindemann}) and~(\ref{eq:Glindemann2}) where the introduction of bound constraints and appropriate regularization helps to introduce desirable qualities into the solutions of phase retrieval problems. 

Several directions for future work also present themselves. Many recent phase retrieval applications from the signal processing literature involve solving large nonlinear least-squares problems like the ones considered in this paper. It would be interesting to extend the ideas presented in this work to those applications. Also, any adaptation of this work for real-world applications would benefit from a more in-depth discussion of regularization parameter selection and appropriate problem-specific regularizers.

Lastly, the codes and methods in this paper are available publicly on Github in the \texttt{BiBox} repository for other researchers looking to use bispectral imaging as part of their work.

%%% Acknowledgements %%%
\section*{Acknowledgements}
\label{sec:Acknowledgements}
We acknowledge funding from the
US National Science Foundation under grant no.\ DMS-1819042 and 1522599.

We would also like to thank Dr. Brandoch Calef at the Maui High Performance Computing Center for helpful advice and guidance.

% Bibliography
\bibliographystyle{IEEEtran}
\bibliography{bispectrum_draft9.bbl}

\begin{thebibliography}{10}
\providecommand{\url}[1]{#1}
\csname url@samestyle\endcsname
\providecommand{\newblock}{\relax}
\providecommand{\bibinfo}[2]{#2}
\providecommand{\BIBentrySTDinterwordspacing}{\spaceskip=0pt\relax}
\providecommand{\BIBentryALTinterwordstretchfactor}{4}
\providecommand{\BIBentryALTinterwordspacing}{\spaceskip=\fontdimen2\font plus
\BIBentryALTinterwordstretchfactor\fontdimen3\font minus
  \fontdimen4\font\relax}
\providecommand{\BIBforeignlanguage}[2]{{%
\expandafter\ifx\csname l@#1\endcsname\relax
\typeout{** WARNING: IEEEtran.bst: No hyphenation pattern has been}%
\typeout{** loaded for the language `#1'. Using the pattern for}%
\typeout{** the default language instead.}%
\else
\language=\csname l@#1\endcsname
\fi
#2}}
\providecommand{\BIBdecl}{\relax}
\BIBdecl

\bibitem{Labeyrie1970}
A.~Labeyrie, ``Attainment of diffraction limited resolution in large telescopes
  by {F}ourier analyzing speckle patterns in star images,'' \emph{Astron.
  Astrophys.}, vol.~6, pp. 85--87, 1970.

\bibitem{Knox1974}
K.~T. Knox and B.~J. Thompson, ``Recovery of images from atmospherically
  degraded short-exposure photographs,'' \emph{Astrophys. J}, vol. 193, pp.
  L45--L48, 1974.

\bibitem{Weigelt1977}
G.~Weigelt, ``Modified speckle interferometry: speckle masking,'' \emph{Opt.
  Commun.}, vol.~21, pp. 55--59, 1977.

\bibitem{Shechtman2015}
Y.~Y.~Shechtman, Y.~C. Eldar, O.~Cohen, H.~N. Chapman, J.~Miao, and M.~Segev,
  ``Phase retrieval with application to optical imaging: a contemporary
  overview,'' \emph{IEEE Signal Processing Magazine}, vol.~32, no.~3, pp.
  87--109, 2015.

\bibitem{Lohmann1983}
A.~W. Lohmann, G.~Weigelt, and B.~Wirnitzer, ``Speckle masking in astronomy:
  triple correlation theory and applications,'' \emph{Applied Optics}, vol.
  22.24, pp. 4028--4037, 1983.

\bibitem{Lohmann1984}
A.~W. Lohmann and B.~Wirnitzer, ``Triple correlations,'' \emph{Proceedings of
  the IEEE}, vol.~72, no.~7, pp. 889--901, 1984.

\bibitem{Negrete1996}
P.~Negrete-Regagnon, ``Practical aspects of image recovery by means of the
  bispectrum,'' \emph{JOSA A}, vol. 13.7, pp. 1557--1576, 1996.

\bibitem{Wirnitzer1985}
B.~Wirnitzer, ``Bispectral analysis at low light levels and astronomical
  speckle masking,'' \emph{JOSA A}, vol.~2, no.~1, pp. 14--21, 1985.

\bibitem{Archer2014}
G.~E. Archer, J.~P. Bos, and M.~C. Roggemann, ``Comparison of bispectrum,
  multiframe blind deconvolution and hybrid bispectrum-multiframe blind
  deconvolution image reconstruction techniques for anisoplanatic, long
  horizontal-path imaging,'' \emph{Optical Engineering}, vol.~53, no.~4, p.
  043109, 2014.

\bibitem{Bos2011}
J.~P. Bos and M.~C. Roggemann, ``Mean squared error performance of
  speckle-imaging using the bispectrum in horizontal imaging applications,'' in
  \emph{Visual Information Processing XX}, vol. 8056.\hskip 1em plus 0.5em
  minus 0.4em\relax International Society for Optics and Photonics, 2011, p.
  805603.

\bibitem{Carrano2002}
C.~J. Carrano, ``Speckle imaging over horizontal paths,'' in
  \emph{High-Resolution Wavefront Control: Methods, Devices, and Applications
  IV}, vol. 4825.\hskip 1em plus 0.5em minus 0.4em\relax International Society
  for Optics and Photonics, 2002, pp. 109--121.

\bibitem{Wen2007}
Z.~Y. Wen, D.~Fraser, A.~Lambert, and H.~D. Li, ``Reconstruction of underwater
  image by bispectrum,'' in \emph{Image Processing, 2007. ICIP 2007. IEEE
  International Conference on}, vol.~3.\hskip 1em plus 0.5em minus 0.4em\relax
  IEEE, 2007, pp. III--545.

\bibitem{Wen2010}
Z.~Y. Wen, A.~Lambert, D.~Fraser, and H.~D. Li, ``Bispectral analysis and
  recovery of images distorted by a moving water surface,'' \emph{Applied
  Optics}, vol.~49, no.~33, pp. 6376--6384, 2010.

\bibitem{Bendory2017}
T.~Bendory, R.~R.~Beinert, and Y.~C. Eldar, ``Fourier phase retrieval:
  uniqueness and algorithms,'' in \emph{Compressed Sensing and its
  Applications}.\hskip 1em plus 0.5em minus 0.4em\relax Springer, 2017, pp.
  55--91.

\bibitem{Bendory2017b}
T.~Bendory, P.~Sidorenko, and Y.~C. Eldar, ``On the uniqueness of frog
  methods,'' \emph{IEEE Signal Processing Letters}, vol.~24, no.~5, pp.
  722--726, 2017.

\bibitem{Bendory2018}
T.~Bendory, Y.~C. Eldar, and N.~Boumal, ``Non-convex phase retrieval from stft
  measurements,'' \emph{IEEE Transactions on Information Theory}, vol.~64,
  no.~1, pp. 467--484, 2018.

\bibitem{Bendory2018b}
T.~Bendory, N.~Boumal, C.~Ma, Z.~Zhao, and A.~Singer, ``Bispectrum inversion
  with application to multireference alignment,'' \emph{IEEE Transactions on
  Signal Processing}, vol.~66, no.~4, pp. 1037--1050, 2018.

\bibitem{Jaganathan2013}
K.~Jaganathan, S.~Oymak, and B.~Hassibi, ``Sparse phase retrieval: convex
  algorithms and limitations,'' in \emph{2013 IEEE International Symposium on
  Information Theory}.\hskip 1em plus 0.5em minus 0.4em\relax IEEE, 2013, pp.
  1022--1026.

\bibitem{Jaganathan2017}
------, ``Sparse phase retrieval: uniqueness guarantees and recovery
  algorithms,'' \emph{IEEE Transactions on Signal Processing}, vol.~65, no.~9,
  pp. 2402--2410, 2017.

\bibitem{Pinilla2018}
S.~Pinilla, J.~Bacca, and H.~Arguello, ``Phase retrieval algorithm via
  nonconvex minimization using a smoothing function,'' \emph{IEEE Transactions
  on Signal Processing}, vol.~66, no.~17, pp. 4574--4584, 2018.

\bibitem{Repetti2014}
A.~Repetti, E.~Chouzenoux, and J.~C. Pesquet, ``A nonconvex regularized
  approach for phase retrieval,'' in \emph{2014 IEEE International Conference
  on Image Processing (ICIP)}.\hskip 1em plus 0.5em minus 0.4em\relax IEEE,
  2014, pp. 1753--1757.

\bibitem{Wang2018}
G.~Wang, L.~Zhang, G.~B. Giannakis, M.~Ak{\c{c}}akaya, and J.~Chen, ``Sparse
  phase retrieval via truncated amplitude flow,'' \emph{IEEE Transactions on
  Signal Processing}, vol.~66, no.~2, pp. 479--491, 2018.

\bibitem{Cho1995}
D.~J. Cho, E.~A. Watson, and G.~M. Morris, ``Application of bispectral speckle
  imaging to near-diffraction-limited imaging in the presence of unknown
  aberrations,'' \emph{Applied Optics}, vol.~34, no.~5, pp. 854--864, 1995.

\bibitem{Kang1991}
M.~G. Kang, K.~T. Lay, and A.~K. Katsaggelos, ``Image restoration algorithms
  based on the bispectrum,'' in \emph{Visual Communications and Image
  Processing'91: Image Processing}, vol. 1606.\hskip 1em plus 0.5em minus
  0.4em\relax International Society for Optics and Photonics, 1991, pp.
  408--419.

\bibitem{Matson1991}
C.~L. Matson, ``Weighted-least-squares phase reconstruction from the
  bispectrum,'' \emph{JOSA A}, vol.~8, no.~12, pp. 1905--1913, 1991.

\bibitem{Meng1990}
J.~Meng, G.~J.~M. Aitken, E.~K. Hege, and J.~S. Morgan, ``Triple-correlation
  subplane reconstruction of photon-address stellar images,'' \emph{JOSA A},
  vol.~7, no.~7, pp. 1243--1250, 1990.

\bibitem{Northcott1988}
M.~J. Northcott, G.~R. Ayers, and J.~C. Dainty, ``Algorithms for image
  reconstruction from photon-limited data using the triple correlation,''
  \emph{JOSA A}, vol.~5, no.~7, pp. 986--992, 1988.

\bibitem{Tyler2004}
D.~Tyler and K.~Schulze, ``Fast phase spectrum estimation using the parallel
  part-bispectrum algorithm,'' \emph{Pub. of the Astron. Soc. of the Pacific},
  vol. 116.815, pp. 65--76, 2004.

\bibitem{Haniff1991}
C.~A. Haniff, ``Least-squares {F}ourier phase estimation from the modulo $2\pi$
  bispectrum phase,'' \emph{JOSA A}, vol. 8.1, pp. 134--140, 1991.

\bibitem{Glindemann1993}
A.~Glindemann and J.~C. Dainty, ``Object fitting to the bispectral phase by
  using least squares,'' \emph{JOSA A}, vol. 10.5, pp. 1056--1063, 1993.

\bibitem{Takajo1988}
H.~Takajo and T.~Takahashi, ``Least-squares phase estimation from the phase
  difference,'' \emph{JOSA A}, vol.~5, no.~3, pp. 416--425, 1988.

\bibitem{Erdem1992}
A.~T. Erdem and M.~I. Sezan, ``Least-squares reconstruction of an image from
  its noisy observations using the bispectrum,'' in \emph{[1992] IEEE Sixth SP
  Workshop on Statistical Signal and Array Processing}.\hskip 1em plus 0.5em
  minus 0.4em\relax IEEE, 1992, pp. 156--159.

\bibitem{Goodman1990}
D.~M. Goodman, T.~W. Lawrence, J.~P. Fitch, and E.~M. Johansson,
  ``Bispectral-based optimization algorithms for speckle imaging,'' in
  \emph{Digital Image Synthesis and Inverse Optics}, vol. 1351.\hskip 1em plus
  0.5em minus 0.4em\relax International Society for Optics and Photonics, 1990,
  pp. 546--561.

\bibitem{Marron1990}
J.~C. Marron, P.~P. Sanches, and R.~C. Sullivan, ``Unwrapping algorithm for
  least-squares phase recovery from the modulo $2\pi$ bispectrum phase,''
  \emph{JOSA A}, vol. 7.1, pp. 14--20, 1990.

\bibitem{Buscher1993}
D.~F. Buscher and C.~A. Haniff, ``Diffraction-limited imaging with partially
  redundant masks: Ii. optical imaging of faint sources,'' \emph{JOSA A},
  vol.~10, no.~9, pp. 1882--1894, 1993.

\bibitem{Nakajima1988}
T.~Nakajima, ``Signal-to-noise ratio of the bispectral analysis of speckle
  interferometry,'' \emph{JOSA A}, vol.~5, no.~9, pp. 1477--1491, 1988.

\bibitem{Roggemann1992}
M.~C. Roggemann and C.~L. Matson, ``Power spectrum and fourier phase spectrum
  estimation by using fully and partially compensating adaptive optics and
  bispectrum postprocessing,'' \emph{JOSA A}, vol.~9, no.~9, pp. 1525--1535,
  1992.

\bibitem{Rudin1992}
L.~I. Rudin, S.~Osher, and E.~Fatemi, ``{Nonlinear total variation based noise
  removal algorithms},'' \emph{Physica D}, vol.~60, no. 1-4, pp. 259--268,
  1992.

\bibitem{NocedalWright1999}
J.~Nocedal and S.~Wright, \emph{Numerical optimization}.\hskip 1em plus 0.5em
  minus 0.4em\relax Berlin: Springer, 1999.

\bibitem{Beck2014}
A.~Beck, \emph{Introduction to Nonlinear Optimization: Theory, Algorithms, and
  Applications with {MATLAB}}.\hskip 1em plus 0.5em minus 0.4em\relax
  Philadelphia, PA: SIAM, 2014.

\bibitem{Haber2014}
E.~Haber, \emph{{Computational Methods in Geophysical Electromagnetics}}.\hskip
  1em plus 0.5em minus 0.4em\relax Philadelphia, PA: SIAM, 2014.

\bibitem{Davis2004}
T.~A. Davis, J.~R. Gilbert, S.~I. Larimore, and E.~G. Ng, ``Algorithm 836:
  Colamd, a column approximate minimum degree ordering algorithm,'' \emph{ACM
  Transactions on Mathematical Software (TOMS)}, vol.~30, no.~3, pp. 377--380,
  2004.

\bibitem{Manteuffel1980}
T.~A. Manteuffel, ``An incomplete factorization technique for positive definite
  linear systems,'' \emph{Mathematics of Computation}, vol.~34, no. 150, pp.
  473--497, 1980.

\bibitem{Saad2003}
Y.~Saad, \emph{Iterative methods for sparse linear systems}.\hskip 1em plus
  0.5em minus 0.4em\relax siam, 2003, vol.~82.

\bibitem{Bjorck1996}
{\AA}.~Bj{\"o}rck, \emph{Numerical methods for least squares problems}.\hskip
  1em plus 0.5em minus 0.4em\relax SIAM, 1996.

\bibitem{golub2013matrix}
G.~H. Golub and C.~F.~V. Loan, \emph{Matrix Computations}, 4th~ed.\hskip 1em
  plus 0.5em minus 0.4em\relax John Hopkins University Press, 2013.

\bibitem{Byrd1995}
R.~H. Byrd, P.~Lu, J.~Nocedal, and C.~Zhu, ``A limited memory algorithm for
  bound constrained optimization,'' \emph{SIAM Journal on Scientific
  Computing}, vol.~16, no.~5, pp. 1190--1208, 1995.

\bibitem{nagy2004iterative}
J.~G. Nagy, K.~M. Palmer, and L.~Perrone, ``Iterative methods for image
  deblurring: {A} {M}atlab object oriented approach,'' \emph{Numerical
  Algorithms}, vol.~36, pp. 73--93, 2004.

\bibitem{RoWe96}
M.~C. Roggemann and B.~M. Welsh, \emph{Imaging through Turbulence}.\hskip 1em
  plus 0.5em minus 0.4em\relax Boca Raton, FL, US: CRC-Press, 1996.

\bibitem{IRtools2018}
S.~Gazzola, P.~C. Hansen, and J.~G. Nagy, ``{IR} {T}ools: {A} {MATLAB} package
  of iterative regularization methods and large-scale test problems,''
  \emph{Numerical Algorithms}, to appear in 2018, doi:
  10.1007/s11075-018-0570-7, see also \url{https://github.com/jnagy1/IRtools}.

\bibitem{MatlabWFS2010}
J.~D. Schmidt, \emph{Numerical Simulation of Optical Wave Propagation with
  Examples in MATLAB}.\hskip 1em plus 0.5em minus 0.4em\relax Bellingham, WA,
  US: SPIE Press, 2010.

\bibitem{Fienup1982}
J.~R. Fienup, ``Phase retrieval algorithms: a comparison,'' \emph{Applied
  Optics}, vol.~21, no.~15, pp. 2758--2769, 1982.

\bibitem{Negrete1996b}
P.~Negrete-Regagnon, ``Phase recovery from the bispectrum aided by the
  error-reduction algorithm,'' \emph{Optics Letters}, vol. 21.4, pp. 275--277,
  1996.

\bibitem{Boyd2004}
S.~Boyd and L.~Vandenberghe, \emph{Convex Optimization}.\hskip 1em plus 0.5em
  minus 0.4em\relax Cambridge University Press, 2004.

\bibitem{Hansen1998}
P.~C. Hansen, \emph{{Rank-deficient and discrete ill-posed problems}}, ser.
  SIAM Monographs on Mathematical Modeling and Computation.\hskip 1em plus
  0.5em minus 0.4em\relax Society for Industrial and Applied Mathematics
  (SIAM), Philadelphia, PA, 1998.

\bibitem{Vogel2002}
C.~R. Vogel, \emph{Computational Methods for Inverse Problems}.\hskip 1em plus
  0.5em minus 0.4em\relax Philadelphia, PA: SIAM, 2002.

\end{thebibliography}

%\begin{IEEEbiography}[{\includegraphics[width=1in,height=1.25in,clip,keepaspectratio]{photo.jpg}}]{James L. Herring}
%\end{IEEEbiography}

%\begin{IEEEbiography}[{\includegraphics[width=1in,height=1.25in,clip,keepaspectratio]{photo.jpg}}]{James Nagy}
%\end{IEEEbiography}

%\begin{IEEEbiography}[{\includegraphics[width=1in,height=1.25in,clip,keepaspectratio]{photo.jpg}}]{Lars Ruthotto}
%\end{IEEEbiography}

\end{document}